# Extending Sicherman Dice to 100-cell Calculation Tables

**Yutaka Nishiyama**
Faculty of Information Management,
Osaka University of Economics,
2, Osumi Higashiyodogawa Osaka, 533-8533, Japan
nishiyama@osaka-ue.ac.jp

**Nozomi Miyanaga**
Faculty of Engineering,
Tamagawa University

**Abstract**

This paper discusses a 100-cell calculation table for addition, $C_{i,j} = A_i + B_j$, $(1 \leq i, j \leq 10)$, addressing the special case where the 100 cells $C_{i,j}$ comprise continuous numbers 0 through 99. We consider what sequences $\{A_i, B_j, (1 \leq i, j \leq 10)\}$ are needed to fulfill these conditions, and the number of possible combinations. We first find seven solutions through factorization of the generating function used for Sicherman dice. We next find seven solutions through geometric means. We demonstrate a correspondence between the seven solutions from factorization and the seven from geometry. Finally, we suggest that the number of solutions $N$ for an $n \times n = n^2$-cell calculation table is $N = (p-2)(p-1) + 1$, where $p$ is the number of divisors of $n$.

**1. The tenth-year miracle**

In 2005, coauthor Nishiyama received an email from Steve Humble—a friend he first met as a visiting fellow at Cambridge—posing an interesting problem. Namely, a pair of six-sided dice with faces labeled 1, 3, 4, 5, 6, 8 and 1, 2, 2, 3, 3, 4 will result in a same distribution of sums as a standard pair of dice with faces labeled 1 through 6. We confirmed that the distributions are indeed the same, and posted the results to the forum of the Japan Mathematical Society (JMS)[1] on 24 May 2005, while Nishiyama was at Cambridge. He later published his findings [6].

From this, coauthor Miyanaga developed the problem into a 100-cell calculation table for addition[2], and published his findings in the same forum on 17 July 2005. He also posed the problem to the students at the university where he is employed, who expressed surprise that factorization could be applied to such a problem. At this point, Miyanaga had applied polynomial factorization to find seven solutions, a method that we will examine in detail below.

Later in 2005, forum posts to the threads "Sums of dice faces" and "100-cell calculations" ceased. Ten years later, in 2015, Miyanaga participated in a symposium hosted by the JMS, where he was surprised to hear a presentation about the same 100-cell addition table and its seven solutions. Even more surprising, this presentation was exactly ten years and one day after his initial post to the forum.

Nishiyama was interested in sums of dice faces (Sicherman dice), but was unfamiliar with 100-cell calculations because they were not used when he was a student. At the time, therefore, he was not particularly interested in Miyanaga's forum post. Upon reviewing that post, however, he realized that this was originally Miyanaga's idea, and urged him to write it up in an article or a research paper. Miyanaga, however, considered the solution using factorization of the generating function applied to Sicherman dice to be of little novelty, and thus not worth an article. An exchange of information and opinions between the two, however, resulted in a surprising development, which we describe later.

---

[1] Est. 2002; website: http://www.sugaku-bunka.org/
[2] A "100-cell calculation table" is a 10 × 10 table with rows and columns commonly labeled with integers 0 through 9 in arbitrary order. An arithmetic operator is specified, and each cell is filled with the result of applying the operator to the row and column labels. This is a popular drill in Japan for students to improve mental calculation speed.

First, we present an overview of the solution demonstrated by Miyanaga.

In the calculation table in Fig. 1, column *A* and row *B* are labeled as follows. This is one solution.

$$A = \{0, 2, 20, 22, 40, 42, 60, 62, 80, 82\}, \qquad B = \{0, 1, 4, 5, 8, 9, 12, 13, 16, 17\}$$

The 100 cells in this table can also be filled in, if preferred (Fig. 2). This can be tedious to do by hand, but a spreadsheet allows rapid calculation through the use of cell references. As one can see, the cells contain integers 0 through 99.

This is one solution to the 100-cell addition table problem, but the more interesting question is how to find the combinations of *A* and *B* labels that give the solution. One method is the factorization used for Sicherman dice, which we describe in detail below.

Nishiyama was playing with the table shown in Fig. 2, wondering whether factorization was the only method of finding solutions, and decided to attempt a geometrical solution. Using the numbers in Fig. 2, he started from 0 and drew a line connecting sequential numbers through 99. As the result in Fig. 3 shows, the path from 0 to 99 showed a pattern. Given that this path represents a solution, finding another path should give another solution. Nishiyama also noticed that the numbers in the leftmost column mirrored the column labels, and the numbers in the uppermost row mirror the row labels (Fig. 4).

We give the details below, but by finding another regular path like that in Fig. 3 and inserting numbers 0 through 99 along that path, taking the second column from the left and second row from the top of Fig. 4 as a new combination of row/column labels, then posing the result as a new problem like in Fig. 1, we should have another specific solution.

| + | 0 | 1 | 4 | 5 | 8 | 9 | 12 | 13 | 16 | 17 |
|---|---|---|---|---|---|---|----|----|----|----|
| 0 |   |   |   |   |   |   |    |    |    |    |
| 2 |   |   |   |   |   |   |    |    |    |    |
| 20 |  |   |   |   |   |   |    |    |    |    |
| 22 |  |   |   |   |   |   |    |    |    |    |
| 40 |  |   |   |   |   |   |    |    |    |    |
| 42 |  |   |   |   |   |   |    |    |    |    |
| 60 |  |   |   |   |   |   |    |    |    |    |
| 62 |  |   |   |   |   |   |    |    |    |    |
| 80 |  |   |   |   |   |   |    |    |    |    |
| 82 |  |   |   |   |   |   |    |    |    |    |

Fig. 1: One solution

| + | 0 | 1 | 4 | 5 | 8 | 9 | 12 | 13 | 16 | 17 |
|---|---|---|---|---|---|---|----|----|----|----|
| 0 | 0 | 1 | 4 | 5 | 8 | 9 | 12 | 13 | 16 | 17 |
| 2 | 2 | 3 | 6 | 7 | 10 | 11 | 14 | 15 | 18 | 19 |
| 20 | 20 | 21 | 24 | 25 | 28 | 29 | 32 | 33 | 36 | 37 |
| 22 | 22 | 23 | 26 | 27 | 30 | 31 | 34 | 35 | 38 | 39 |
| 40 | 40 | 41 | 44 | 45 | 48 | 49 | 52 | 53 | 56 | 57 |
| 42 | 42 | 43 | 46 | 47 | 50 | 51 | 54 | 55 | 58 | 59 |
| 60 | 60 | 61 | 64 | 65 | 68 | 69 | 72 | 73 | 76 | 77 |
| 62 | 62 | 63 | 66 | 67 | 70 | 71 | 74 | 75 | 78 | 79 |
| 80 | 80 | 81 | 84 | 85 | 88 | 89 | 92 | 93 | 96 | 97 |
| 82 | 82 | 83 | 86 | 87 | 90 | 91 | 94 | 95 | 98 | 99 |

Fig. 2: Values filled in

Fig. 3: Cells 0 through 99 connected

Fig. 4: Leftmost and topmost columns

Nishiyama focused on pathing through association with the knight's tour problem and general solutions to odd × odd magic squares. Indeed, since the knight's tour involves moving a knight so that it visits each of the 64 squares on a chess board exactly once, a solution to that problem is also a solution for an 8 × 8 addition table. Figure 5a shows one such solution, with the cell marked "1" the starting point and 2, 3, …, 64 the path. Figure 5b shows the first 25 moves in that path as connected lines.

A magic square problem involves finding an arrangement of numbers such that rows, columns, and diagonals have the same sums, and general solutions are known for odd × odd and even × even magic squares. Figure 6 shows a solution for an odd × odd (5 × 5) magic square, with 1 placed in the center of the bottom row and the path proceeding diagonally, down and to the left, moving to cells in corresponding locations when spilling over the grid. When a diagonal cell already containing a number is encountered, we move straight up instead. This allows a sequential path from 1 through 25, following a pathing rule. Since the addition table problem for 0 through 99 too is a problem related to finding paths, Nishiyama was convinced there must be some geometric method to finding solutions.

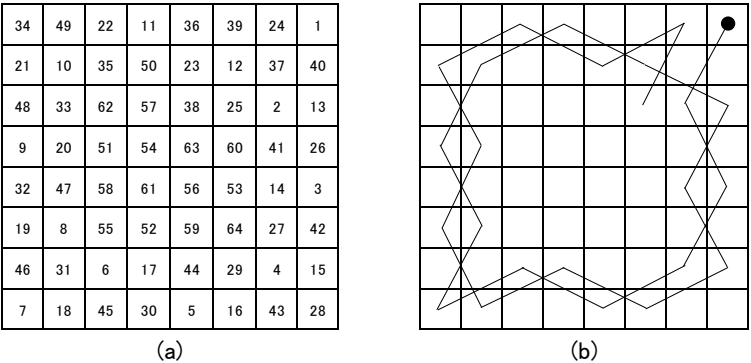

Fig. 5: A knight's tour solution [4]

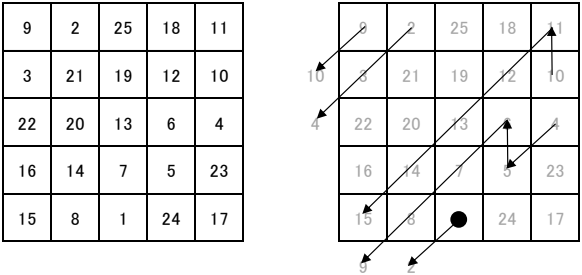

Fig. 6: A general solution to the magic squares problem (odd × odd)

**2. Sicherman dice**

This section gives a simple explanation of Sicherman dice. These are a special kind of dice developed by Col. George Sicherman, which have the same distribution of face sums (Fig. 7b) as do a canonical pair of dice (Fig. 7a). The faces on a pair of Sicherman dice are labeled 1, 3, 4, 5, 6, 8 and 1, 2, 2, 3, 3, 4. Martin Gardner discussed these dice in a 1978 *Scientific American* article [1].

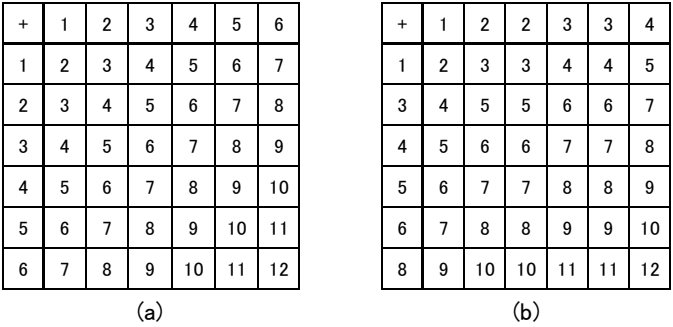

Fig. 7: Distributions of sums are the same for (a) standard and (b) Sicherman dice

Gallian et al. applied factorization of a generating function for Sicherman dice to show that they are the only possible labeling that give the same distribution as standard dice [2][3]. Specifically, in the polynomial factorization below, the second line represents standard dice and the third line represents Sicherman dice. See Ref. [6] for details.

$$x^2 + 2x^3 + 3x^4 + 4x^5 + 5x^6 + 6x^7 + 5x^8 + 4x^9 + 3x^{10} + 2x^{11} + x^{12}$$
$$= (x + x^2 + x^3 + x^4 + x^5 + x^6)(x + x^2 + x^3 + x^4 + x^5 + x^6)$$
$$= (x + x^3 + x^4 + x^5 + x^6 + x^8)(x + 2x^2 + 2x^3 + x^4)$$

There's a bit more to add regarding Sicherman dice. I had my article [6] translated into English [7], allowing it to be read by people all over the world. In particular, in January 2012 I was surprised to receive an email from George Sicherman himself. I was suspicious at first that it was truly him, but I learned that he had a web page, and was able to confirm his identity. He had written to thank me for contributing to the popularity of Sicherman dice, and to comment on a line I had written at the end of my paper:

> These days the dice are known as 'crazy dice', or alternatively they take the name of their discoverer, i.e., 'Sicherman dice'. Sicherman sells these dice as merchandise, although it seems that they are not used in actual casinos. [7]

As it turns out Sicherman himself does not sell these dice, though a company called GameStation and several others do, using his name without permission (Fig. 8). I wrote back to notify him that I had updated this line in both my Japanese and English papers to reflect the true state of affairs.

Sicherman dice have retained their popularity, and remain topical even today. As recently as 2014 the *New York Times* wrote an article about them [9], in which they mentioned my English paper.

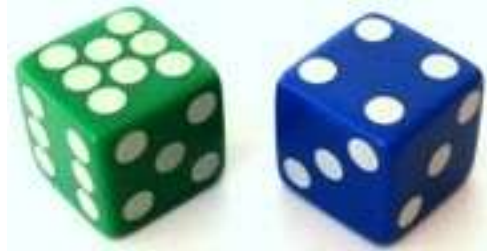

Fig. 8: Commercially available Sicherman dice (GameStation)

**3. Seven solutions from factorization**

This section describes the method of solution through factorization that Miyanaga applied to the 100-cell calculation table for addition.

We want to add $A$ and $B$ to obtain $C$, with $C$ being $\{0, 1, 2, \ldots, 99\}$. The problem is to find $A$ and $B$, so we can represent $A$, $B$, $C$ as $a(x), b(x), c(x)$, as per the following equation.

$$a(x) \times b(x) = c(x)$$
$$c(x) = 1 + x + x^2 + x^3 + \cdots + x^{99}$$

Here, $c(x)$ is called a generating function. The polynomial is difficult to factorize as-is, so we multiply both sides by $(1 - x)$.

$$c(x) \times (1 - x) = (1 + x + x^2 + x^3 + \cdots + x^{99}) \times (1 - x) = 1 - x^{100}$$

We then try factorization for $1 - x^{100}$.

$$1 - x^{100} = (1 - x)(1 + x)(1 + x^2)(1 + x + x^2 + x^3 + x^4)(1 - x + x^2 - x^3 + x^4)$$
$$(1 - x^2 + x^4 - x^6 + x^8)(1 + x^5 + x^{10} + x^{15} + x^{20})(1 - x^5 + x^{10} - x^{15} + x^{20})$$
$$(1 - x^{10} + x^{20} - x^{30} + x^{40})$$

Since we multiplied this by $(1 - x)$, we can divide by $(1 - x)$ to obtain $c(x)$.

$$c(x) = (1 + x)(1 + x^2)(1 + x + x^2 + x^3 + x^4)(1 - x + x^2 - x^3 + x^4)(1 - x^2 + x^4 - x^6 + x^8)$$
$$(1 + x^5 + x^{10} + x^{15} + x^{20})(1 - x^5 + x^{10} - x^{15} + x^{20})(1 - x^{10} + x^{20} - x^{30} + x^{40})$$

Thus, $c(x)$ has been factorized into eight terms. Sorting these terms into $a(x)$ and $b(x)$ groups and expanding them gives sequences for $A$ and $B$ labels for a 100-cell addition table.

Doing this, Miyanaga found seven solutions. Nishiyama independently searched for solutions using factorization, and comparison of the results found them to be in perfect agreement with Miyanaga's. Dividing the eight terms into two groups and confirming that each group is the sum of ten terms without omission would require an enormous amount of calculation. To ensure that there are no omissions, let us substitute $x = 1$ into the polynomial and try a method similar to that used for Sicherman dice [6].

Say that $c(x)$ is factored into eight functions $p(x)$ to $w(x)$. We substitute $x = 1$ into these equations and examine the resulting values.

$$c(x) = 1 + x + x^2 + x^3 + \cdots + x^{99} = p(x)q(x)r(x)s(x)t(x)u(x)v(x)w(x)$$
$$\begin{aligned} p(x) &= 1 + x, & p(1) &= 2 \\ q(x) &= 1 + x^2, & q(1) &= 2 \\ r(x) &= 1 + x + x^2 + x^3 + x^4, & r(1) &= 5 \\ s(x) &= 1 + x^5 + x^{10} + x^{15} + x^{20}, & s(1) &= 5 \\ t(x) &= 1 - x + x^2 - x^3 + x^4, & t(1) &= 1 \\ u(x) &= 1 - x^2 + x^4 - x^6 + x^8, & u(1) &= 1 \\ v(x) &= 1 - x^5 + x^{10} - x^{15} + x^{20}, & v(1) &= 1 \\ w(x) &= 1 - x^{10} + x^{20} - x^{30} + x^{40}, & w(1) &= 1 \end{aligned}$$

Assume that polynomials $\{a(x), b(x)\}$ fulfill the following conditions:
- Condition 1: 10 terms in $a(x)$, each with coefficient 1
- Condition 2: 10 terms in $b(x)$, each with coefficient 1
- Condition 3: $a(x) \times b(x) = c(x)$

From conditions 1 and 2, it is necessary that

$$c(1) = 100, \quad a(1) = 10, \quad b(1) = 10.$$

Also,

$$p(1) = q(1) = 2, \quad r(1) = s(1) = 5$$
$$p(1) \times r(1) = p(1) \times s(1) = 2 \times 5 = 10,$$

so from Condition 3 the combination of $\{p(x), q(x), r(x), s(x)\}$ divided among $\{a(x), b(x)\}$ can be

$$\{p(x) \times r(x), q(x) \times s(x)\} \text{ or } \{p(x) \times s(x), q(x) \times r(x)\}.$$

Meanwhile,

$$t(1) = u(1) = v(1) = w(1) = 1,$$

so it is possible that $\{t(x), u(x), v(x), w(x)\}$ are assigned to either $a(x)$ or $b(x)$. From this, there are $2 \times 2^4 = 32$ possible assignments of $\{p(x), q(x), r(x), s(x), t(x), u(x), v(x), w(x)\}$ to $\{a(x), b(x)\}$. After assignment, we must expand each of the terms and carefully ensure that requirements 1–3 hold. This is possible using paper and pencil, but that is a lot of work and prone to error, so we can also use mathematical processing software such as Wolfram Alpha[3].

Below are the seven solutions found using this method.

*Solution 1*
$a(x) = (1 + x^2)(1 - x^2 + x^4 - x^6 + x^8)(1 + x^5 + x^{10} + x^{15} + x^{20})(1 - x^5 + x^{10} - x^{15} + x^{20})$
$(1 - x^{10} + x^{20} - x^{30} + x^{40})$
$\quad = (1 + x^2)(1 - x^2 + x^4 - x^6 + x^8)(1 + x^{10} + x^{20} + x^{30} + x^{40})(1 - x^{10} + x^{20} - x^{30} + x^{40})$
$\quad = (1 + x^{10})(1 + x^{20} + x^{40} + x^{60} + x^{80})$
$\quad = 1 + x^{10} + x^{20} + x^{30} + x^{40} + x^{50} + x^{60} + x^{70} + x^{80} + x^{90}$
$b(x) = (1 + x)(1 + x + x^2 + x^3 + x^4)(1 - x + x^2 - x^3 + x^4)$
$\quad = (1 + x)(1 + x^2 + x^4 + x^6 + x^8)$
$\quad = 1 + x + x^2 + x^3 + x^4 + x^5 + x^6 + x^7 + x^8 + x^9$

---
[3] http://www.wolframalpha.com/

*Solution 2*
$$a(x) = (1+x)(1-x+x^2-x^3+x^4)(1+x^5+x^{10}+x^{15}+x^{20})(1-x^5+x^{10}-x^{15}+x^{20})$$
$$= 1+x^5+x^{10}+x^{15}+x^{20}+x^{25}+x^{30}+x^{35}+x^{40}+x^{45}$$
$$b(x) = (1+x^2)(1+x+x^2+x^3+x^4)(1-x^2+x^4-x^6+x^8)(1-x^{10}+x^{20}-x^{30}+x^{40})$$
$$= 1+x+x^2+x^3+x^4+x^{50}+x^{51}+x^{52}+x^{53}+x^{54}$$

*Solution 3*
$$a(x) = (1+x^2)(1-x^2+x^4-x^6+x^8)(1-x^{10}+x^{20}-x^{30}+x^{40})(1+x^5+x^{10}+x^{15}+x^{20})$$
$$= 1+x^5+x^{10}+x^{15}+x^{20}+x^{50}+x^{55}+x^{60}+x^{65}+x^{70}$$
$$b(x) = (1+x)(1+x+x^2+x^3+x^4)(1-x+x^2-x^3+x^4)(1-x^5+x^{10}-x^{15}+x^{20})$$
$$= 1+x+x^2+x^3+x^4+x^{25}+x^{26}+x^{27}+x^{28}+x^{29}$$

*Solution 4*
$$a(x) = (1+x)(1-x+x^2-x^3+x^4)(1+x^5+x^{10}+x^{15}+x^{20})(1-x^5+x^{10}-x^{15}+x^{20})$$
$$(1-x^{10}+x^{20}-x^{30}+x^{40})$$
$$= 1+x^5+x^{20}+x^{25}+x^{40}+x^{45}+x^{60}+x^{65}+x^{80}+x^{85}$$
$$b(x) = (1+x^2)(1-x^2+x^4-x^6+x^8)(1+x+x^2+x^3+x^4)$$
$$= 1+x+x^2+x^3+x^4+x^{10}+x^{11}+x^{12}+x^{13}+x^{14}$$

*Solution 5*
$$a(x) = (1+x)(1+x^5+x^{10}+x^{15}+x^{20})(1-x^5+x^{10}-x^{15}+x^{20})(1-x^{10}+x^{20}-x^{30}+x^{40})$$
$$= 1+x+x^{20}+x^{21}+x^{40}+x^{41}+x^{60}+x^{61}+x^{80}+x^{81}$$
$$b(x) = (1+x^2)(1+x+x^2+x^3+x^4)(1-x+x^2-x^3+x^4)(1-x^2+x^4-x^6+x^8)$$
$$= 1+x^2+x^4+x^6+x^8+x^{10}+x^{12}+x^{14}+x^{16}+x^{18}$$

*Solution 6*
$$a(x) = (1+x^2)(1-x^2+x^4-x^6+x^8)(1-x^{10}+x^{20}-x^{30}+x^{40})(1+x+x^2+x^3+x^4)$$
$$(1-x+x^2-x^3+x^4)$$
$$= 1+x^2+x^4+x^6+x^8+x^{50}+x^{52}+x^{54}+x^{56}+x^{58}$$
$$b(x) = (1+x)(1+x^5+x^{10}+x^{15}+x^{20})(1-x^5+x^{10}-x^{15}+x^{20})$$
$$= 1+x+x^{10}+x^{11}+x^{20}+x^{21}+x^{30}+x^{31}+x^{40}+x^{41}$$

*Solution 7*
$$a(x) = (1+x^2)(1+x^5+x^{10}+x^{15}+x^{20})(1-x^5+x^{10}-x^{15}+x^{20})(1-x^{10}+x^{20}-x^{30}+x^{40})$$
$$= 1+x^2+x^{20}+x^{22}+x^{40}+x^{42}+x^{60}+x^{62}+x^{80}+x^{82}$$
$$b(x) = (1+x)(1+x+x^2+x^3+x^4)(1-x+x^2-x^3+x^4)(1-x^2+x^4-x^6+x^8)$$
$$= 1+x+x^4+x^5+x^8+x^9+x^{12}+x^{13}+x^{16}+x^{17}$$

In the above, we have divided the eight terms between $a(x)$ and $b(x)$ and verified that they are the sum of ten terms, but let us use geometric progressions to verify that $a(x) \times b(x) = c(x)$.

*Solution 1*

$a(x)$ has an initial term of 1, common ratio $x^{10}$, and is a sum of a geometric progression of ten terms, while $b(x)$ has an initial term of 1, common ratio $x$, and is a sum of a geometric progression of ten terms, so
$$a(x) = 1+x^{10}+x^{20}+x^{30}+x^{40}+x^{50}+x^{60}+x^{70}+x^{80}+x^{90}$$
$$= \frac{1 \times (1-(x^{10})^{10})}{1-x^{10}} = \frac{1-x^{100}}{1-x^{10}}$$
$$b(x) = 1+x+x^2+x^3+x^4+x^5+x^6+x^7+x^8+x^9$$
$$= \frac{1 \times (1-x^{10})}{1-x} = \frac{1-x^{10}}{1-x}.$$
From this,
$$a(x) \times b(x) = \frac{1-x^{100}}{1-x^{10}} \times \frac{1-x^{10}}{1-x} = \frac{1-x^{100}}{1-x} = c(x).$$

Similarly,

*Solution 2*

$$a(x) = 1 + x^5 + x^{10} + x^{15} + x^{20} + x^{25} + x^{30} + x^{35} + x^{40} + x^{45}$$
$$= \frac{1 \times (1-(x^5)^{10})}{1-x^5} = \frac{1-x^{50}}{1-x^5}$$
$$b(x) = 1 + x + x^2 + x^3 + x^4 + x^{50} + x^{51} + x^{52} + x^{53} + x^{54}$$
$$= (1 + x + x^2 + x^3 + x^4)(1+x^{50})$$
$$= \frac{1 \times (1-x^5)}{1-x} \times (1+x^{50}) = \frac{1-x^5}{1-x} \times \frac{1-x^{100}}{1-x^{50}}$$
$$a(x) \times b(x) = \frac{1-x^{50}}{1-x^5} \times \frac{1-x^5}{1-x} \times \frac{1-x^{100}}{1-x^{50}} = \frac{1-x^{100}}{1-x} = c(x)$$

*Solution 3*

$$a(x) = 1 + x^5 + x^{10} + x^{15} + x^{20} + x^{50} + x^{55} + x^{60} + x^{65} + x^{70}$$
$$= (1 + x^5 + x^{10} + x^{15} + x^{20})(1+x^{50})$$
$$= \frac{1 \times (1-(x^5)^5)}{1-x^5} \times (1+x^{50}) = \frac{1-x^{25}}{1-x^5} \times \frac{1-x^{100}}{1-x^{50}}$$
$$b(x) = 1 + x + x^2 + x^3 + x^4 + x^{25} + x^{26} + x^{27} + x^{28} + x^{29}$$
$$= (1 + x + x^2 + x^3 + x^4)(1+x^{25})$$
$$= \frac{1 \times (1-x^5)}{1-x} \times (1+x^{25}) = \frac{1-x^5}{1-x} \times \frac{1-x^{50}}{1-x^{25}}$$
$$a(x) \times b(x) = \frac{1-x^{25}}{1-x^5} \times \frac{1-x^{100}}{1-x^{50}} \times \frac{1-x^5}{1-x} \times \frac{1-x^{50}}{1-x^{25}} = \frac{1-x^{100}}{1-x} = c(x)$$

*Solution 4*
$$a(x) = 1 + x^5 + x^{20} + x^{25} + x^{40} + x^{45} + x^{60} + x^{65} + x^{80} + x^{85}$$
$$= (1 + x^{20} + x^{40} + x^{60} + x^{80})(1+x^5)$$
$$= \frac{1 \times (1-(x^{20})^5)}{1-x^{20}} \times (1+x^5) = \frac{1-x^{100}}{1-x^{20}} \times \frac{1-x^{10}}{1-x^5}$$
$$b(x) = 1 + x + x^2 + x^3 + x^4 + x^{10} + x^{11} + x^{12} + x^{13} + x^{14}$$
$$= (1 + x + x^2 + x^3 + x^4)(1+x^{10})$$
$$= \frac{1 \times (1-x^5)}{1-x} \times (1+x^{10}) = \frac{1-x^5}{1-x} \times \frac{1-x^{20}}{1-x^{10}}$$
$$a(x) \times b(x) = \frac{1-x^{100}}{1-x^{20}} \times \frac{1-x^{10}}{1-x^5} \times \frac{1-x^5}{1-x} \times \frac{1-x^{20}}{1-x^{10}} = \frac{1-x^{100}}{1-x} = c(x)$$

*Solution 5*
$$a(x) = 1 + x + x^{20} + x^{21} + x^{40} + x^{41} + x^{60} + x^{61} + x^{80} + x^{81}$$
$$= (1 + x^{20} + x^{40} + x^{60} + x^{80})(1+x)$$
$$= \frac{1 \times (1-(x^{20})^5)}{1-x^{20}} \times (1+x) = \frac{1-x^{100}}{1-x^{20}} \times (1+x)$$
$$b(x) = 1 + x^2 + x^4 + x^6 + x^8 + x^{10} + x^{12} + x^{14} + x^{16} + x^{18}$$
$$= \frac{1 \times (1-(x^2)^{10})}{1-x^2} = \frac{1-x^{20}}{1-x^2}$$
$$a(x) \times b(x) = \frac{1-x^{100}}{1-x^{20}} \times (1+x) \times \frac{1-x^{20}}{1-x^2} = \frac{1-x^{100}}{1-x} = c(x)$$

*Solution 6*
$$a(x) = 1 + x^2 + x^4 + x^6 + x^8 + x^{50} + x^{52} + x^{54} + x^{56} + x^{58}$$
$$= (1 + x^2 + x^4 + x^6 + x^8)(1+x^{50})$$
$$= \frac{1 \times (1-(x^2)^5)}{1-x^2} \times (1+x^{50}) = \frac{1-x^{10}}{1-x^2} \times \frac{1-x^{100}}{1-x^{50}}$$
$$b(x) = 1 + x + x^{10} + x^{11} + x^{20} + x^{21} + x^{30} + x^{31} + x^{40} + x^{41}$$
$$= (1 + x^{10} + x^{20} + x^{30} + x^{40})(1+x)$$

$$= \frac{1 \times (1-(x^{10})^5)}{1-x^{10}} \times (1+x) = \frac{1-x^{50}}{1-x^{10}} \times (1+x)$$
$$a(x) \times b(x) = \frac{1-x^{10}}{1-x^2} \times \frac{1-x^{100}}{1-x^{50}} \times \frac{1-x^{50}}{1-x^{10}} \times (1+x) = \frac{1-x^{100}}{1-x} = c(x)$$

*Solution 7*
$a(x) = 1 + x^2 + x^{20} + x^{22} + x^{40} + x^{42} + x^{60} + x^{62} + x^{80} + x^{82}$
$\phantom{a(x)} = (1 + x^{20} + x^{40} + x^{60} + x^{80})(1 + x^2)$
$$= \frac{1 \times (1-(x^{20})^5)}{1-x^{20}} \times (1+x^2) = \frac{1-x^{100}}{1-x^{20}} \times (1+x^2)$$
$b(x) = 1 + x + x^4 + x^5 + x^8 + x^9 + x^{12} + x^{13} + x^{16} + x^{17}$
$\phantom{b(x)} = (1 + x^4 + x^8 + x^{12} + x^{16})(1 + x)$
$$= \frac{1 \times (1-(x^4)^5)}{1-x^4} \times (1+x) = \frac{1-x^{20}}{1-x^4} \times (1+x)$$
$$a(x) \times b(x) = \frac{1-x^{100}}{1-x^{20}} \times (1+x^2) \times \frac{1-x^{20}}{1-x^4} \times (1+x) = \frac{1-x^{100}}{1-x} = c(x)$$

A recent publication presents a problem similar to the 100-cell addition table (with $n = 4$) [8]. The proof uses a generating function to show that there are three possible factorizations and is attributed to "Hiroshi Yoshikawa, in a 2006 seminar for third-year students at the Tokyo Institute of Technology, Dept. of Information Science," so it seems likely that the method of applying Sicherman dice factorization provided a hint to the solution.

$$\begin{aligned}1 + x + x^2 + \cdots + x^{15} &= (1+x)(1+x^2)(1+x^4)(1+x^8) \\ &= (1+x+x^2+x^3)(1+x^4+x^8+x^{12}) \\ &= (1+x+x^4+x^5)(1+x^2+x^8+x^{10}) \\ &= (1+x+x^8+x^9)(1+x^2+x^4+x^6)\end{aligned}$$

**4. Cyclotomic polynomials**

We would like to examine the factorization of $1 - x^{100}$ a little more closely. The formula

$$a^2 - b^2 = (a+b)(a-b),$$

which is frequently encountered in high school, allows factorization into three terms, as

$$1 - x^{100} = (1+x^{50})(1-x^{50}) = (1+x^{50})(1+x^{25})(1-x^{25}).$$

These terms can be factorized as follows:

$$\begin{aligned}(1+x^{50}) &= (1+x^{10})(1-x^{10}+x^{20}-x^{30}+x^{40}), \\ (1+x^{25}) &= (1+x^5)(1-x^5+x^{10}-x^{15}+x^{20}), \\ (1-x^{25}) &= (1-x^5)(1+x^5+x^{10}+x^{15}+x^{20}).\end{aligned}$$

Furthermore,

$$\begin{aligned}(1+x^{10}) &= (1+x^2)(1-x^2+x^4-x^6+x^8), \\ (1+x^5) &= (1+x)(1-x+x^2-x^3+x^4), \\ (1-x^5) &= (1-x)(1+x+x^2+x^3+x^4),\end{aligned}$$

so substitution allows decomposition into nine terms:

$$\begin{aligned}1-x^{100} = (1-x)(1+x)(1+x^2)(1+x+x^2+x^3+x^4)(1-x+x^2-x^3+x^4) \\ (1-x^2+x^4-x^6+x^8)(1+x^5+x^{10}+x^{15}+x^{20})(1-x^5+x^{10}-x^{15}+x^{20}) \\ (1-x^{10}+x^{20}-x^{30}+x^{40}).\end{aligned}$$

We multiplied $c(x)$ by $(1-x)$ to make factorization easier, so dividing it back out gives eight terms for

$c(x)$.

We can use methods such as this to factorize even high-degree polynomials, but the special case of polynomials in the form $1 - x^n$ have been discussed in the context of cyclotomic polynomials, so let us take a look at their relation.

There is a relation between $x^{100} - 1 = 0$ and the factorization of $1 - x^{100}$. As Nishiyama has written about in the past, there is also a relation between $x^{17} - 1 = 0$ and Gauss' method for constructing a 17-sided regular polyhedron [5]. $x^n - 1 = 0$ is called the cyclotomic equation, and states that there are $n$ complex roots including $x = 1$ on a unit circle.

Knowing the equation for a cyclotomic polynomial, it is relatively easy to perform high-dimensional factorizations. Cyclotomic polynomials are related to roots of unity, and are defined as $\Phi_n(x)$ as

$$\Phi_n(x) = \frac{x^n - 1}{\prod_{d|n, d \neq n} \Phi_d(x)}.$$

This is a polynomial with integer coefficients, and is an irreducible polynomial over the rationals. A polynomial $x^n - 1$ can be irreducibly decomposed as a product of cyclotomic polynomials as

$$x^n - 1 = \prod_{d|n} \Phi_d(x),$$

where $d$ is a divisor of $n$.

Some actual calculations of cyclotomic polynomials follow:

$\Phi_1 = \phantom{(x^2 - 1)/\Phi_1} = x - 1$

$\Phi_2 = (x^2 - 1)/\Phi_1 \phantom{xx} = x + 1$

$\Phi_3 = (x^3 - 1)/\Phi_1 \phantom{xx} = x^2 + x + 1$

$\Phi_4 = (x^4 - 1)/\Phi_1 \Phi_2 \phantom{x} = x^2 + 1$

$\Phi_5 = (x^5 - 1)/\Phi_1 \phantom{xx} = x^4 + x^3 + x^2 + x + 1$

$\Phi_6 = (x^6 - 1)/\Phi_1 \Phi_2 \Phi_3 = x^2 - x + 1$

$\Phi_7 = (x^7 - 1)/\Phi_1 \phantom{xx} = x^6 + x^5 + x^4 + x^3 + x^2 + x + 1$

$\Phi_8 = (x^8 - 1)/\Phi_1 \Phi_2 \Phi_4 = x^4 + 1$

$\Phi_9 = (x^9 - 1)/\Phi_1 \Phi_3 \phantom{x} = x^6 + x^3 + 1$

$\Phi_{10} = (x^{10} - 1)/\Phi_1 \Phi_2 \Phi_5 = x^4 - x^3 + x^2 - x + 1$

We require a factorization of $x^{100} - 1$, so following the second equation above we have

$$x^{100} - 1 = \Phi_1 \Phi_2 \Phi_4 \Phi_5 \Phi_{10} \Phi_{20} \Phi_{25} \Phi_{50} \Phi_{100}.$$

Here, $n = 100$ and 100 has nine divisors $d$ (1, 2, 4, 5, 10, 20, 25, 50, 100), each a product of irreducible polynomial $\Phi_d$. We can perform substitutions such as

$\Phi_{20} = (x^{20} - 1)/\Phi_1 \Phi_2 \Phi_4 \Phi_5 \Phi_{10} \phantom{xxxxxxx} = x^8 - x^6 + x^4 - x^2 + 1$
$\Phi_{25} = (x^{25} - 1)/\Phi_1 \Phi_5 \phantom{xx} = x^{20} + x^{15} + x^{10} + x^5 + 1$
$\Phi_{50} = (x^{50} - 1)/\Phi_1 \Phi_2 \Phi_5 \Phi_{10} \Phi_{25} \phantom{xxxxxx} = x^{20} - x^{15} + x^{10} - x^5 + 1$
$\Phi_{100} = (x^{100} - 1)/\Phi_1 \Phi_2 \Phi_4 \Phi_5 \Phi_{10} \Phi_{20} \Phi_{25} \Phi_{50} = x^{40} - x^{30} + x^{20} - x^{10} + 1,$

giving

$$x^{100} - 1 = (x - 1)(x + 1)(x^2 + 1)(x^4 + x^3 + x^2 + x + 1)(x^4 - x^3 + x^2 - x + 1)$$

$$(x^8 - x^6 + x^4 - x^2 + 1)(x^{20} + x^{15} + x^{10} + x^5 + 1)(x^{20} - x^{15} + x^{10} - x^5 + 1)$$
$$(x^{40} - x^{30} + x^{20} - x^{10} + 1),$$

so while the order of presentation of the $c(x) \times (1 - x)$ above slightly differs, we have obtained a similar factorization.

**5. Geometric solution**

This section describes an alternate method of solving the 100-cell addition problem, one employing a geometric approach.

Let $c(x) = a(x) \times b(x)$ be a $10 \times 10 = 100$-cell grid with cell 0 in the upper-left corner and 99 in the lower-right corner. We also have rules for pathing from cell 0 to 99. We also apply these rules between blocks, as described below.

・Rule 1: Move from left to right.
・Rule 2: Move from above to below.

Let us consider ways of partitioning these $10 \times 10 = 100$ cells into even blocks. 10 has four divisors, 1, 2, 5, and 10. Ignoring 10 as a divisor gives us three ways of partitioning the blocks, and since these can be applied either horizontally or vertically we have $3 \times 3 = 9$ ways of partitioning in total. However, after applying the above two rules and excluding repeated patterns (like those in Fig. 9), we are left with 7 partitions.

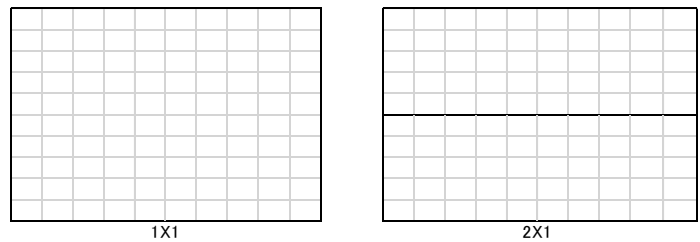

Fig. 9: Excluded same solutions

Figure 10 shows the partitions. The cells in the topmost partitioning have not been split at all, those in the middle row have been split in half, and those in the bottom row have been split into five. Those in the middle and bottom row are in three patterns, those with no splitting, those split in half, and those split into five equal parts. In all, there are $1 + 3 + 3 = 7$ patterns.

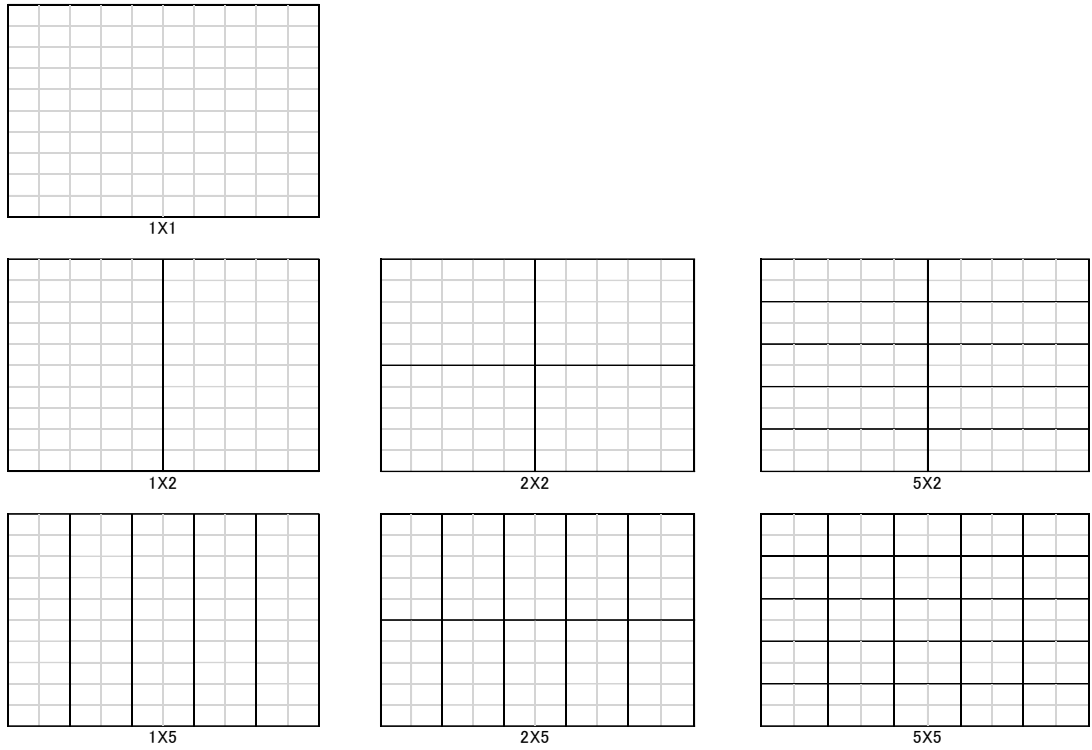
Fig. 10: Seven patterns for division

Figure 11 shows the results of filling in numbers 0 through 99 along the paths resulting from our rules. The paths in the central figures are not only within blocks, but between blocks as well. Right-side figures show numbers 0 through 99 filled in, with leftmost columns in each representing $a(x)$ and topmost rows representing $b(x)$. This is a solution to the 100-cell addition table problem.

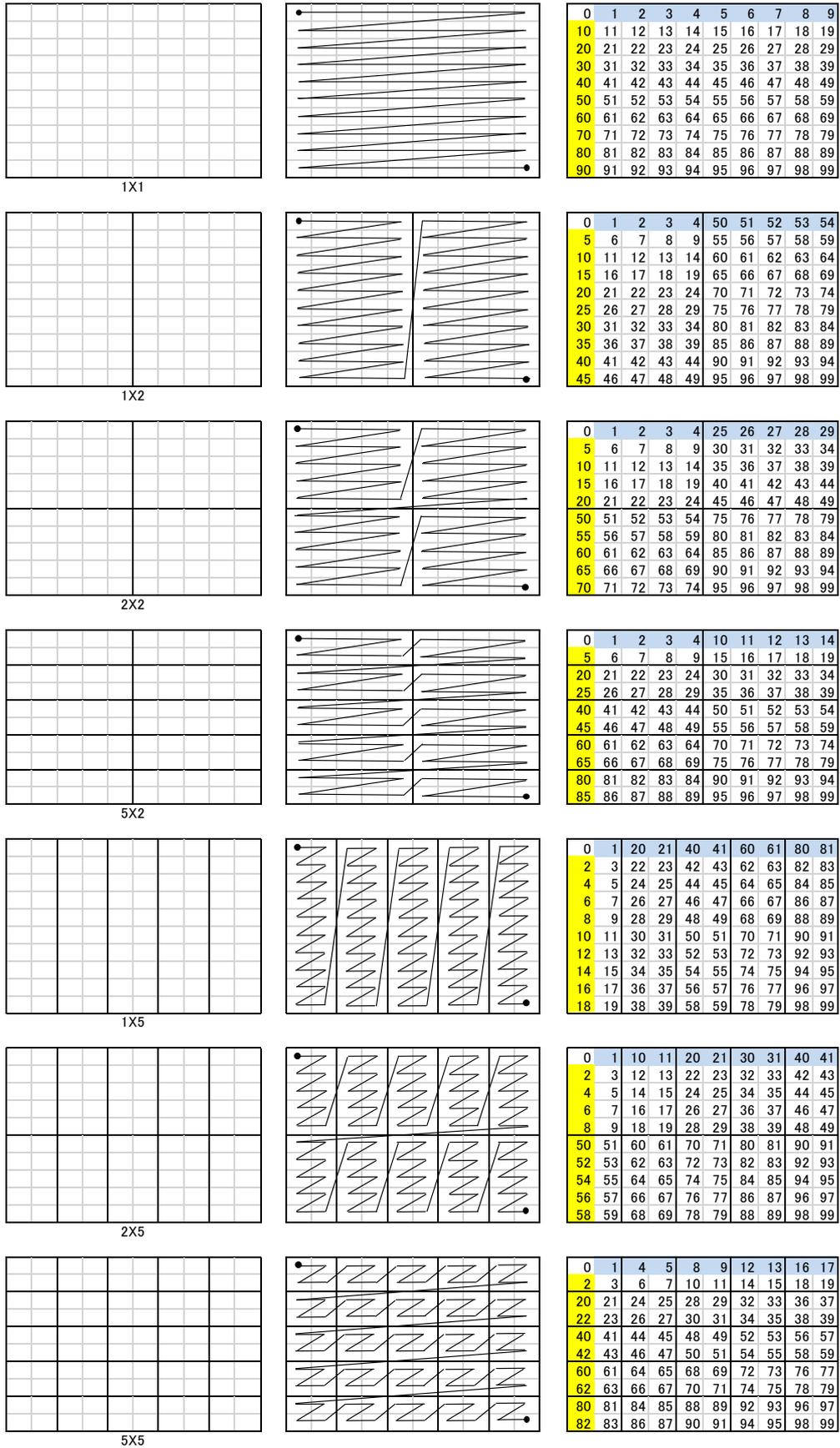

Fig. 11: Seven solutions from geometry

Letting leftmost sequences be *A* and uppermost rows be *B*, we get the following seven solutions. We can associate these geometric solutions with the solutions from factorization $a(x), b(x)$.

*Solution 1*
$$A = \{0, 10, 20, 30, 40, 50, 60, 70, 80, 90\}, B = \{0, 1, 2, 3, 4, 5, 6, 7, 8, 9\}$$
$$a(x) = 1 + x^{10} + x^{20} + x^{30} + x^{40} + x^{50} + x^{60} + x^{70} + x^{80} + x^{90}$$
$$b(x) = 1 + x + x^2 + x^3 + x^4 + x^5 + x^6 + x^7 + x^8 + x^9$$

*Solution 2*
$$A = \{0, 5, 10, 15, 20, 25, 30, 35, 40, 45\}, B = \{0, 1, 2, 3, 4, 50, 51, 52, 53, 54\}$$
$$a(x) = 1 + x^5 + x^{10} + x^{15} + x^{20} + x^{25} + x^{30} + x^{35} + x^{40} + x^{45}$$
$$b(x) = 1 + x + x^2 + x^3 + x^4 + x^{50} + x^{51} + x^{52} + x^{53} + x^{54}$$

*Solution 3*
$$A = \{0, 5, 10, 15, 20, 50, 55, 60, 65, 70\}, B = \{0, 1, 2, 3, 4, 25, 26, 27, 28, 29\}$$
$$a(x) = 1 + x^5 + x^{10} + x^{15} + x^{20} + x^{50} + x^{55} + x^{60} + x^{65} + x^{70}$$
$$b(x) = 1 + x + x^2 + x^3 + x^4 + x^{25} + x^{26} + x^{27} + x^{28} + x^{29}$$

*Solution 4*
$$A = \{0, 5, 20, 25, 40, 45, 60, 65, 80, 85\}, B = \{0, 1, 2, 3, 4, 10, 11, 12, 13, 14\}$$
$$a(x) = 1 + x^5 + x^{20} + x^{25} + x^{40} + x^{45} + x^{60} + x^{65} + x^{80} + x^{85}$$
$$b(x) = 1 + x + x^2 + x^3 + x^4 + x^{10} + x^{11} + x^{12} + x^{13} + x^{14}$$

*Solution 5*
$$A = \{0, 1, 20, 21, 40, 41, 60, 61, 80, 81\}, B = \{0, 2, 4, 6, 8, 10, 12, 14, 16, 18\}$$
$$a(x) = 1 + x + x^{20} + x^{21} + x^{40} + x^{41} + x^{60} + x^{61} + x^{80} + x^{81}$$
$$b(x) = 1 + x^2 + x^4 + x^6 + x^8 + x^{10} + x^{12} + x^{14} + x^{16} + x^{18}$$

*Solution 6*
$$A = \{0, 2, 4, 6, 8, 50, 52, 54, 56, 58\}, B = \{0, 1, 10, 11, 20, 21, 30, 31, 40, 41\}$$
$$a(x) = 1 + x^2 + x^4 + x^6 + x^8 + x^{50} + x^{52} + x^{54} + x^{56} + x^{58}$$
$$b(x) = 1 + x + x^{10} + x^{11} + x^{20} + x^{21} + x^{30} + x^{31} + x^{40} + x^{41}$$

*Solution 7*
$$A = \{0, 2, 20, 22, 40, 42, 60, 62, 80, 82\}, B = \{0, 1, 4, 5, 8, 9, 12, 13, 16, 17\}$$
$$a(x) = 1 + x^2 + x^{20} + x^{22} + x^{40} + x^{42} + x^{60} + x^{62} + x^{80} + x^{82}$$
$$b(x) = 1 + x + x^4 + x^5 + x^8 + x^9 + x^{12} + x^{13} + x^{16} + x^{17}$$

## 6. Generalization as $N = (p-2)(p-1) + 1$

We have seen that there are seven solutions for a $10 \times 10 = 100$-cell calculation table. So how many solutions will there be for an $n \times n = n^2$-cell table? This section considers how we might generalize *n*.

First, we examine the geometric solutions for $2 \times 2 = 4$ to $9 \times 9 = 81$-cell tables.

Partitioning into evenly sized blocks was the key to finding the number of solutions for 100-cell tables, so we show these divisions in Fig. 12.

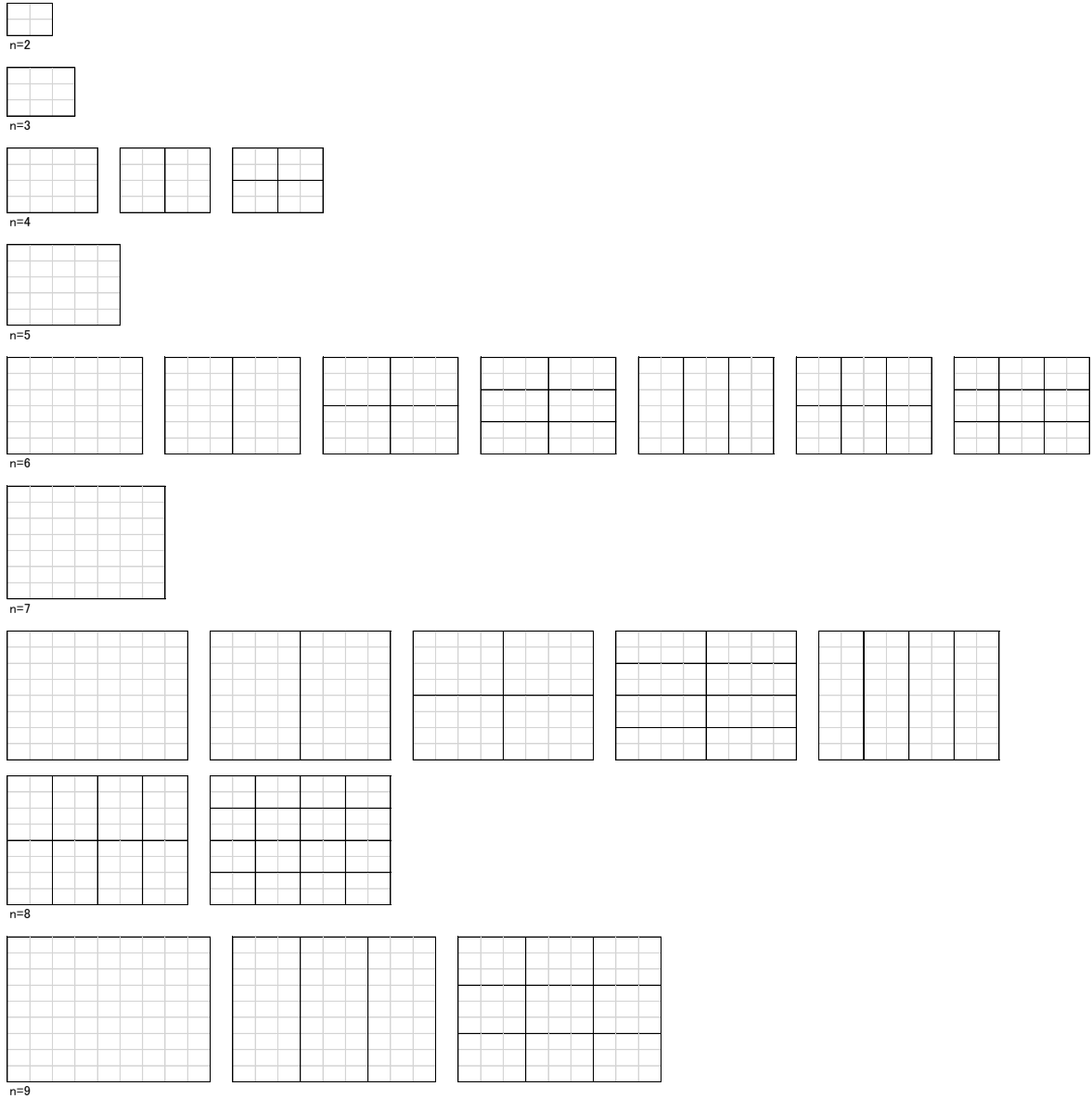

Fig. 12: Evenly dividing $n \times n = n^2$-cell tables into blocks ($2 \leq n \leq 9$)

| $n$ | $n \times n = n^2$ | Divisors | Num. divisors ($p$) | Num. solutions ($N$) |
|---|---|---|---|---|
| 2 | 4 | 1, 2 | 2 | 1 |
| 3 | 9 | 1, 3 | 2 | 1 |
| 4 | 16 | 1, 2, 4 | 3 | 3 |
| 5 | 25 | 1, 5 | 2 | 1 |
| 6 | 36 | 1, 2, 3, 6 | 4 | 7 |
| 7 | 49 | 1, 7 | 2 | 1 |
| 8 | 64 | 1, 2, 4, 8 | 4 | 7 |
| 9 | 81 | 1, 3, 9 | 3 | 3 |
| 10 | 100 | 1, 2, 5, 10 | 4 | 7 |

Table 1. Number of divisors ($p$) and number of solutions ($N$)

First, we note that if $n$ is prime ($n = 2, 3, 5, 7$) the blocks cannot be evenly divided, so there is only one solution. In the case of $n = 4, 9$ there are three possible partitions, so there are three solutions, and in the case of $n = 6, 8$ there are seven possible partitions, so there are seven solutions. Solutions using factorization confirmed these results.

Divisors of *n* are related to the number of possible even partitions. Table 1 lists the divisors of *n* and the number of divisors *p*.

In the case of an $n \times n = n^2$-cell table, letting the number of solutions be *N*, the relation between *N* and *n* is as follows:

$$\begin{aligned} N &= 1 & (n = 2, 3, 5, 7), \\ N &= 3 & (n = 4, 9), \\ N &= 7 & (n = 6, 8, 10). \end{aligned}$$

We can rewrite these *N* values as follows:

$$N = 1, \quad N = 3 = 1 + 2, \quad N = 7 = 1 + 3 + 3.$$

Considering the case of $n = 12$ for a $12 \times 12 = 144$-cell table, there are six divisors of 12 (1, 2, 3, 4, 6, 12), so we can expect the number of solutions to be

$$N = 21 = 1 + 5 + 5 + 5 + 5, \quad (n = 12).$$

We omit the details here, but again the geometric and factorization methods agreed that there are indeed 21 solutions.

This led to the following definition and prediction.

**Definition**

Let *A* and *B* be sequences of length *n*, $A = \{a_1, a_2, a_3, \cdots, a_n\}, B = \{b_1, b_2, b_3, \cdots, b_n\}$, and let $A_i = a_i$ and $B_j = b_j$. Let $C = \{C_{i,j}\}$ be calculated through addition of *A* and *B* as

$$C_{i,j} = A_i + B_j, \quad (1 \leq i, j \leq n.)$$

Even if the $n \times n = n^2$ elements of *C* are randomly ordered, the numbers 0 through $n^2 - 1$ can be arranged as sums in the $n^2$-cell calculation table without duplication.

**Prediction**

Letting *p* be the number of divisors of *n*, the number of solutions *N* for an $n \times n = n^2$-cell calculation table is

$$N = (p - 2)(p - 1) + 1.$$

Ex: For $n = 12$, there are 6 divisors (1, 2, 3, 4, 6, 12), so $p = 6$ and $N = 4 \times 5 + 1 = 21$.
For $n = 16$, there are 5 divisors (1, 2, 4, 8, 16), so $p = 5$ and $N = 3 \times 4 + 1 = 13$.

The above was Nishiyama's prediction, but when trying to develop a proof Miyanaga found a counterexample. For example, in the case of $n = 8$ a solution like that in Fig. 13 is possible even with the same division pattern. In other words, when there is a nested structure of blocks within blocks. There are likely other solutions, which we leave for more detailed consideration in a future paper.

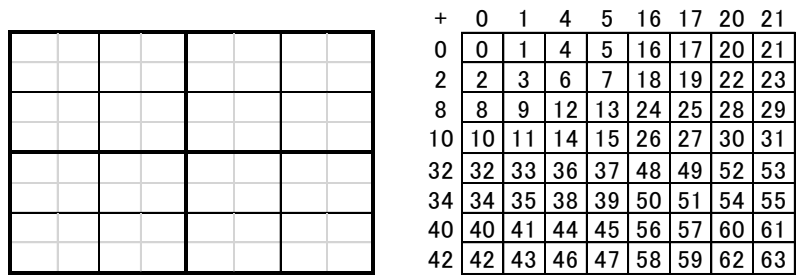

Fig. 13: Alternate solution in a nested structure ($n = 8$)